\documentclass[12pt]{article}
\usepackage{amsmath}
\newcommand{\be}{\begin{equation}}
\newcommand{\ee}{\end{equation}}
\newcommand{\ba}{\begin{eqnarray}}
\newcommand{\ea}{\end{eqnarray}}
\newcommand{\ban}{\begin{eqnarray*}}
\newcommand{\ean}{\end{eqnarray*}}

 \newcommand{\qed}{\hspace*{\fill}\rule{3mm}{3mm}\quad}
\newcommand{\Pf}{\noindent  {\em Proof.} }

\newcommand{\sect}[1]{\section{#1}  \setcounter{equation}{0}}

\newtheorem{lem}{Lemma}[section]
 
\begin{document}
\newtheorem{defn}[lem]{Definition}
\newtheorem{theo}[lem]{Theorem}
\newtheorem{cor}[lem]{Corollary}
\newtheorem{prop}[lem]{Proposition}
\newtheorem{rk}[lem]{Remark}
\newtheorem{ex}[lem]{Example}
\newtheorem{note}[lem]{Note}
\newtheorem{conj}[lem]{Conjecture}

\title{Sobolev Inequalities, Riesz Transforms  and the Ricci Flow} 
\author{Rugang Ye \\ {\small Department  of Mathematics} \\
{\small University of California, Santa  Barbara}}
\date{August 25, 2007}
\maketitle

\noindent 1. Introduction \\
2. From $W^{1,p}$ for  lower $p$ to $W^{1,p}$ for higher $p$\\
3. Nonlocal Sobolev inequalities in terms of the $(1,p)$-Bessel norm \\
4. $W^{2,p}$ Sobolev inequalities \\
5. Estimates of the Riesz transform and $W^{1,p}$ Sobolev inequalities\\

\sect{Introduction} 

\noindent {\bf Part A} \\

In [Y1], [Y2], [Y3] and [Y4], logarithmic Sobolev inequalties along the Ricc flow were 
obtained. As a consequence, $W^{1,p}$ Sobolev inequalities with $p=2$ along the Ricci flow were derived. 
Let $M$ be a closed manifold of dimension $n \ge 2$. Let $g=g(t)$ be a smooth solution of the Ricci flow on $M \times [0, T)$ for some (finite or infinite) $T>0$.
In the case $\lambda_0(g_0) \ge 0$, where $\lambda_0(g_0)$ denotes the first eigenvalue of 
the operator $-\Delta+\frac{R}{4}$ of the initial metric $g_0=g(0)$, the Sobolev inequality takes the following form
(in the case $n \ge 3$)
\ba \label{sobolev1}
\left(\int_M |u|^{\frac{2n}{n-2}} dvol\right)^{\frac{n-2}{n}} \le A \int_M (|\nabla u|^2 
+\frac{R}{4}u^2)dvol,
\ea
where the constant $A$ depends on the initial metric in terms of rudimentary 
geometric data. If the condition $\lambda_0(g_0)\ge 0$ is not assumed, then
the Sobolev inequality takes the form (again in the case $n\ge 3$)
\ba \label{sobolev2}
\left(\int_M |u|^{\frac{2n}{n-2}} dvol\right)^{\frac{n-2}{n}} \le A \int_M (|\nabla u|^2 
+\frac{R}{4}u^2)dvol+B\int_M u^2 dvol,
\ea
where the constants $A$ and $B$ depend on a finite upper bound for $T$ and the initial metric 
in terms of rudimentary geometric data.

As is well-known, the case $p=2$ of the $W^{1,p}(M)$ Sobolev inequalities is the most important for applications to analysis and geometry. However, it is of high 
interest to understand the situation $1<p<2$ and $2<p<n$, both from the point of view of 
a deeper understanding of the theory and the point of view of further applications.  
In this paper, we derive $W^{1,p}$ and $W^{2,p}$ Sobolev inequalities for general $p$ along the Ricci flow in several different ways. We'll take a general point of view, and study the general problem of deriving 
further Sobolev inequalities from a given Sobolev inequality.   In particular, we'll include 
noncompact manifolds and manifolds with boundary, which require additional care.  
Our first result is the following one.  \\

\noindent {\bf Theorem A1} {\it  Let $(M, g)$ be a  Riemannian 
manifold of dimension $n \ge 2$, with or without boundary. (It is not assumed to be compact or complete.)   Assume for some $1\le p_0<n$
\ba \label{p0}
\left(\int_M |u|^{\frac{p_0n}{n-p_0}} dvol\right)^{\frac{n}{n-p_0}} \le
A \int_M |\nabla u|^{p_0} dvol+\frac{B}{vol_g(M)^{\frac{p_0}{n}}} \int_M |u|^{p_0} dvol.
\ea
Then we have for all $p_0<p<n$
\ba \label{pp}
\left(\int_M |u|^{\frac{np}{n-p}}dvol\right)^{\frac{n-p}{n}} \le 
C_1 \int_M |\nabla u|^p dvol 
+\frac{C_2}{vol_g(M)^{\frac{p}{n}}} \int_M |u|^p dvol,
\ea
where the constants $C_1=C_1(n,p_0, p,A,B)$ and $C_2=C_2(n,p_0, p,A,B)$ depend only on 
$n,p_0, p ,A$ and $B$. Their dependence on $p$ is in terms of an upper bound for 
$\frac{1}{n-p}$.  (If $vol_g(M)=\infty$, then it is understood that everything 
involving $B$ is nonpresent.)} \\

This theorem is proved by an induction scheme based on the H\"{o}lder inequality. 
The principle that Sobolev inequalties of lower $p$ lead to Sobolev inequalities of 
higher $p$ is known. For example, it is well-known that the Sobolev inequality 
\ba
\|u\|_{\frac{np}{n-p}} \le C \|\nabla u\|_p 
\ea
for $u \in W^{1,p}_0({\bf R}^n)$, $1<p<n$, can be derived from the case $p=1$, see e.g. [GT].
However, the result in Theorem A1 is new, and the 
proof is more involved. Combining this theorem with the results in [Y1], [Y2] and 
[Y3] we then obtain $W^{1,p}$ Sobolev inequalitites along the Ricci flow 
for $p>2$ in various sitautions.  To keep this paper streamlined, we only state 
the results in the situation of [Y1]. The results in the 
situations of [Y2] and [Y3] are similar and obvious.   \\

\noindent {\bf Theorem A2} {\it  Assume that $\lambda_0(g_0) > 0$. Let $2<p<n$.  
Then there holds for each $t \in [0, T)$ and 
all $u \in W^{1,p}(M)$  
\ba \label{sobolevD}
\left( \int_M |u|^{\frac{np}{n-p}} dvol \right)^{\frac{n-p}{n}} \le 
A\left[(\max R^++1)vol(M)^{\frac{2}{n}}\right]^{\frac{m(p) p}{2}} \int_M (|\nabla u|^p+|u|^p) dvol,
\nonumber \\
\ea
where all geometric quantities are associated with $g(t)$,  except the constant $A$, 
which can be bounded from above in terms of the dimension $n$, a nonpositive lower bound for $R_{g_0}$, a positive lower bound for $vol_{g_0}(M)$, an upper bound for $C_S(M,g_0)$,
a positive lower bound for $\lambda_0(g_0)$, and an upper bound for $\frac{1}{n-p}$. The quantity 
$(\max R^++1)vol(M)^{\frac{2}{n}}$ is at time $t$ and the 
number $m(p)$ is defined in the proof of this theorem below. } \\

\noindent {\bf Theorem A3} {\it   Assume $T<\infty$. 
Let $2<p<n$.  
Then there holds for each $t \in [0, T)$ and 
all $u \in W^{1,2}(M)$  
\ba \label{sobolevD*}
\left( \int_M |u|^{\frac{np}{n-p}} dvol \right)^{\frac{n-p}{n}} \le 
A\left[ 1+(\max R^++1)vol(M)^{\frac{2}{n}}\right]^{\frac{m(p)p}{2}} \int_M (|\nabla u|^p+|u|^p)dvol,
\nonumber \\
\ea
where all geometric quantities are associated with $g(t)$,  except the constant 
$A$, which  can be bounded from above in terms of the dimension $n$, 
a nonpositive lower bound for $R_{g_0}$, a positive lower bound for $vol_{g_0}(M)$, an upper bound for $C_S(M,g_0)$,
an upper bound for $T$, and an upper bound for $\frac{1}{n-p}$. 
The quantity 
$(\max R^++1)vol(M)^{\frac{2}{n}}$ is at time $t$ and the 
number $m(p)$ is the same as in Theorem A2. } \\

\noindent {\bf Theorem A4} {\it Let $n=3$ and $g=g(t)$ be a Ricci flow with surgeries as constructed 
in [P2] on its maximal time interval $[0, T_{max})$, with suitably chosen surgery parameters.  Let $g_0=g(0)$. Let $m(t)$ denote the number of 
surgeries which are performed up to the time $t\in (0, T_{max})$. Let $2<p<3$.  Then there holds at  each $t \in [0, T_{max})$
\ba \label{sobolevJ*}
\left( \int_M |u|^{\frac{3p}{p-3}} dvol \right)^{\frac{p-3}{3}} \le 
A(t)\left[ 1+(\max R^++1)vol(M)^{\frac{2}{3}}\right]^{\frac{m_p p}{2}} \int_M (|\nabla u|^p+|u|^p)dvol
\nonumber \\
\ea
for all $u \in W^{1,2}(M)$, where $A(t)$ is bounded from above in terms of a nonpositive lower bound for $R_{g_0}$, a positive lower bound for $vol_{g_0}(M)$, an upper bound for $C_S(M,g_0)$, 
an upper bound for $t$, and an upper bound for $\frac{1}{3-p}$.  The quantity 
$(\max R^++1)vol(M)^{\frac{2}{n}}$ is at time $t$ and the 
number $m(p)$ is the same as in Theorem A2. 

If $\lambda_0(g_0)>0$, then $A(t)$ can be bounded from above in terms of  a nonpositive lower bound for $R_{g_0}$, a positive lower bound for $vol_{g_0}(M)$, an upper bound for $C_S(M,g_0)$,
a positive lower bound for $\lambda_0(g_0)$, an upper bound for $m(t)$, and 
an upper bound for $\frac{1}{3-p}$.} 
\\

The remaining results for the Ricci flow in this paper extend to the Ricci flow with surgeries in the 
same fashion as this last theorem. We omit the statements of these extensions. \\

\noindent {\bf Part B} \\

First we present a  result on nonlocal Sobolev inequalities which is analogous to 
[Y1, Theorem C.5], but is formulated in terms the  canonical $(1,p)$-Bessel norm for Sobolev 
functions. \\

\noindent {\bf Definition} Let $(M, g)$ be a metrically complete Riemannian manifold, with or without boundary. 
Let $1<p<\infty$. Let the Bessel-Sobolev space $L_B^{1,p}(M)$ be the completion of 
$C^{\infty}_c(M)$ with respect to the norm $\|(-\Delta+1)^{\frac{1}{2}}\|_p$. 
(See Section 3 for the constrcution of the operator $(-\Delta+1)^{\frac{1}{2}}$. )  
We shall say that $g$ is a {\it $p$-Bessel metric}, and $(M,g)$ is a {\it $p$-Bessel (Riemannian) manifold}, 
if $L^{1,p}_B(M)$ is equivalent to $W^{1,p}(M)$, i.e. 
\ba
c_1 \|u\|_{1,p} \le \|(-\Delta+1)^{\frac{1}{2}}u\|_p 
\le c_2 \|u\|_{1,p}
\ea
for all $u \in C^{\infty}_{c}(M)$ and some positive constants $c_1$ and $c_2$, where 
$\|u\|_{1,p}=\|u\|_p+\|\nabla u\|_p$ is the $W^{1,p}$ norm of $u$.   

Assume that $(M,g)$ is $p$-Bessel. We define the $(1,p)$-Bessel norm for $f \in W^{1,p}(M)$ 
to be 
$\|f\|_{B, 1, p} =\|(-\Delta+1)^{\frac{1}{2}} f\|_p.$
 \\

The operators $(-\Delta+1)^{\alpha}$ are called Bessel potentials, which is the reason for the 
above terminologies involving ``Bessel".  Note that every metrically complete $(M ,g)$ is 
$2$-Bessel because of the identity 
\ba
<(-\Delta+1)^{\frac{1}{2}}u, (-\Delta+1)^{\frac{1}{2}}u>_2
=\int_M (|\nabla u|^2+u^2)dvol
\ea
for $u\in C^{\infty}_c(M)$ (and then also for $u \in W^{1,2}(M)$).  By the arguments in [Y1, Appendix C] (see also [St]), $(M,g)$ is $p$-Bessel for each $1<p<\infty$  if $M$ is compact.  
\\

\noindent {\bf Theorem B1} {\it  Let $(M,g)$ be a metrically complete manifold with or without boundary of dimension 
$n \ge 2$. Let $1<\mu<\infty$.  Assume  the Sobolev inequality 
\ba \label{SB1}
\left(\int_M |u|^{\frac{2\mu}{\mu-2}}dvol\right)^{\frac{\mu-2}{\mu}} \le A\int_M |\nabla u|^{2}dvol
+B \int_M |u|^{2}dvol
\ea 
for all $u \in W^{1,2}(M)$.
Let $1<p<\mu$.  Then there  
holds
\ba \label{SB1-p}
\|(-\Delta+1)^{-\frac{1}{2}}u\|_{\frac{\mu p}{\mu-p}} \le C(\mu, A, B, p) \|u\|_p
\ea
for all $u \in L^p(M)$, where the constant $C(\mu, A, B, p)$ 
can be bounded from above in terms of upper bounds for $A$, $B$, $\mu$,  
$\frac{1}{\mu-p}$ and $\frac{1}{p-1}$. Consequently, there holds for a given $1<p<\mu$
\ba \label{SB1-pp}
\|u\|_{\frac{\mu p}{\mu -p}} \le C(\mu, A, B, p) \|u\|_{B,1,p}
\ea
for all $u \in W^{1,p}(M)$, provided that $(M, g)$ is $p$-Bessel. } \\

Combining this theorem with the results in [Y1],[Y2] and [Y3] we obtain 
nonlocal Sobolev inequalitities along the Ricci flow which are analogous to 
Theorem C.6 and Theorem C.7 in [Y1], but are formulated in terms of the 
canonical $(1,p)$-Bessel norm. Again we only state the results in the situation of 
[Y1]. As before, let $g=g(t)$ be a smooth solution of the Ricci flow on $M \times [0, T)$ 
for a closed manifold $M$ of dimension $n \ge 3$ and some $0<T \le \infty$. 
\\

\noindent {\bf Theorem B2} {\it Assume that $\lambda_0(g_0)>0$.  Let $1<p<n$.   There is a  positive constant $C$ depending only on 
the dimension $n$, a positive lower bound for $\lambda_0(g_0)$, a positive lower bound for $vol_{g_0}(M)$, an upper bound for $C_S(M,g_0)$, an upper bound for 
$\frac{1}{p-1}$, and an upper bound for $\frac{1}{n-p}$, such that for each $t \in [0, T)$ and 
all $u \in W^{1,p}(M)$ there holds 
\ba  \label{SB2-p}
\|u\|_{\frac{np}{n-p}} \le C (1+R_{max}^+)^{\frac{1}{2}}\|u\|_{B,1,p}.
\ea
} \\

\noindent {\bf Theorem B3} {\it Assume $T<\infty$ and $1<p<n$.  There is a positive constant $C$ depending only on the dimension $n$, 
a nonpositive lower bound for $R_{g_0}$, a positive lower bound for $vol_{g_0}(M)$, an upper bound for $C_S(M,g_0)$,
an upper bound for $T$, an upper bound for $\frac{1}{p-1}$, and an upper bound 
for $\frac{1}{n-p}$,  such that for each $t \in [0, T)$ and 
all $u \in W^{1,p}(M)$ there holds 
\ba \label{SB3-p}
\|u\|_{\frac{np}{n-p}} \le C (1+ R_{max}^+)^{\frac{1}{2}} \|u\|_{B,1,p}.
\ea
} \\

The inequality (\ref{SB1-p}) in Theorem B-1 is a special case of the following more general result. \\

\noindent {\bf Theorem B4} {\it Let $(M,g)$ be a metrically complete manifold, possibly with boundary.
Let $\Psi \in L^{\infty}(M)$ and $\mu>1$. Assume that $\Psi\ge 0$ and  
the Sobolev inequality 
\ba \label{SB4}
\left(\int_M |u|^{\frac{2\mu}{\mu-2}} dvol \right)^{\frac{\mu-2}{\mu}} 
\le A\int_M (|\nabla u|^2+\Psi u^2)dvol
\ea
for some $A>0$.  Set $H=-\Delta+\Psi$.  Let $1<p<\mu$.  Then  
there holds 
\ba \label{SB4-p}
\|H^{-\frac{1}{2}}u\|_{\frac{\mu p}{\mu-p}} \le C(A, \mu, p) \|u\|_p  
\ea
for all $u \in L^p(M)$, where the positive constant 
 $C(\mu, c, p)$ can be bounded from above in terms of upper bounds for $A$, $\mu$,  
$\frac{1}{\mu-p}$ and $\frac{1}{p-1}$.  Consequently, there holds
for each $1<p<\infty$ 
\ba \label{SB4-pp}
\|u\|_{\frac{\mu p}{\mu-p}} \le C(A, \mu, p) \|H^{\frac{1}{2}}u\|_p
\ea
for all $u \in W^{1,p}(M)$, provided that $(M, g)$ is compact.
}\\

\noindent {\bf Part C} \\

Next we have the following consequence of Theorem B4.
\\

\noindent {\bf Theorem C1} {\it Let $(M, g)$ be a compact  Riemannian manifold of dimension $n\ge 2$,
with or without boundary. Assume $\Psi \in L^{\infty}(M)$ and $\mu>1$. 
Assume $\Psi\ge 0$ and the Sobolev inequality
\ba \label{SC1}
\left(\int_M |u|^{\frac{2\mu}{\mu-2}} dvol\right)^{\frac{\mu-2}{\mu}} \le
A \int_M (|\nabla u|^2+\Psi u^2) dvol.  \ea
 Let $1<p<\frac{\mu}{2}$.  Then  
there holds 
\ba \label{SC1-pp}
\|u\|_{\frac{\mu p}{\mu-2p}} \le C(\mu,A, p) \|\Delta u+\Psi u\|_p  
\ea
for all $u \in W^{2,p}(M)$, where the constant 
 $C(\mu,A, p)$ can be bounded from above in terms of upper bounds for $\mu$,
 $A$, 
$\frac{1}{\mu-p}$ and $\frac{1}{p-1}$. } \\

Combinging this theorem with the results in [Y1],[Y2] and [Y3] 
we then obtain $W^{2,p}$ Sobolev inequalities along the Ricci flow.
Again, we only state the results in the situation of [Y1]. Let 
$g=g(t)$ be a smooth solution of the Ricci flow on $M \times [0, T)$ 
for a closed manifold of dimension $n \ge 3$ and some $0<T \le \infty$, 
with a given initial metric $g_0$. \\

\noindent {\bf Theorem C2} {\it Assume that $R_{g_0} \ge 0$ and $\lambda_0(g_0)>0$ (thus $R_{g_0}$ is somewhere 
positive).  Let $1<p<\frac{n}{2}$.   There is a  positive constant $C$ depending only on 
the dimension $n$, a positive lower bound for $\lambda_0(g_0)$, a positive lower bound for $vol_{g_0}(M)$, an upper bound for $C_S(M,g_0)$, an upper bound for 
$\frac{1}{p-1}$, and an upper bound for $\frac{1}{n-2p}$, such that for each $t \in [0, T)$ there holds 
\ba  \label{SB2-p}
\|u\|_{\frac{np}{n-2p}} \le C \|\Delta u-\frac{R}{4}u\|_{p}
\ea
for all $u \in W^{2, p}(M)$.
} \\

\noindent {\bf Theorem C3} {\it Assume $T<\infty$ and $1<p<n$.  There is a positive constant $C$ depending only on the dimension $n$, 
a nonpositive lower bound for $R_{g_0}$, a positive lower bound for $vol_{g_0}(M)$, an upper bound for $C_S(M,g_0)$,
an upper bound for $T$, an upper bound for $\frac{1}{p-1}$, and an upper bound 
for $\frac{1}{n-2p}$,  such that for each $t \in [0, T)$ there holds 
\ba \label{SB3-p}
\|u\|_{\frac{np}{n-2p}} \le C \|\Delta u-(\frac{R}{4}-\frac{\min R^-_{g_0}}{4}+1)u\|_{p}
\ea
for all $u \in W^{2, p}(M)$.
} \\

\noindent {\bf Part D} \\

Now we address the issue of converting the nonlocal $W^{1,p}$ Sobolev inequalitites in Part B into
conventional $W^{1,p}$ Sobolev inequalities.  Let $H$ denote the operator $-\Delta+\Psi$. Obviously, 
the desired convertion requires an estimate of the following kind 
\ba \label{conversion}
\|H^{\frac{1}{2}}u\|_p \le C \|u\|_{1,p}
\ea
for all $u \in W^{1,p}(M)$. 
Assume that $M$ is compact. Since $H$ is a pseudo-differential operator of order 1 [Se], the inequality  
(\ref{conversion})  holds true for some $C$, as mentioned before for the special case 
$\Psi=1$.    But 
the constant $C$ obtained this way depends on 
$M$ and the metric $g$ in rather complicated ways. Our purpose  is 
to obtain a constant $C$ which has clear and rudimentary geometric dependences. 
For this purpose, the general theory 
of pseudo-differential operators does not seem to give any information.

The issue at hand can be understood in terms of the Riesz trasform  of $H$, which is defined to be 
${R}_H=\nabla H^{-\frac{1}{2}}$.  An $L^p$ inequality for the Riesz transform 
\ba \label{riesz}
\|{R}_Hu\|_p \le c \|u\|_p
\ea
for all $u \in L^p(M)$
means the same as 
\ba \label{riesz1}
\|\nabla u\|_p \le c \|H^{\frac{1}{2}}u\|_p
\ea
for all $u \in L^p(M)$.
On the other hand, by duality, the inequality 
(\ref{riesz1}) implies (\ref{conversion}) for the dual exponent 
under suitable conditions on $\Psi$. (In the special case $\Psi=0$, it leads to 
$\|H^{\frac{1}{2}}u\|_q \le c\|\nabla u\|_q$
for the dual exponent $q$.) 
In general, the Riesz transform ${R}_L$ of a nonnegative symmetric elliptic operator 
$L$ (of second order) is defined in the same way as ${R}_H$.    
A fundamental problem in harmonic analysis and potential theory is to obtain $L^p$ boundedness for Riesz tranforms ${R}_L$, or 
inequalities $\|\nabla u\|_p \le C \|L^{\frac{1}{2}}u\|_p$ and 
$\|L^{\frac{1}{2}}u\|_p \le C\|u\|_{1,p}$.  
 From a geometric point of view, 
$L^p$-boundedness alone is not enough. It is crucial to obtain geometric estimates for 
the constants.    

Based on the $L^p$ estimates for Riesz transforms due to D.~Bakry [B] we can convert 
the nonlocal $W^{1,p}$ Sobolev inequalities in Theorem B1 to conventional 
$W^{1,p}$ Sobolev inequalities which depend on a lower bound for the Ricci curvature. 
The situation of Theorem B4 with a general $\Psi$ is more complicated, our corresponding result will 
be presented elsewhere.\\

\noindent {\bf Theorem D1} {\it  Let $(M,g)$ be a complete manifold (without boundary) of dimension 
$n \ge 2$. Let $1<\mu<\infty$.  Assume  the Sobolev inequality 
\ba \label{SD1}
\left(\int_M |u|^{\frac{2\mu}{\mu-2}}dvol\right)^{\frac{\mu-2}{\mu}} \le A\int_M |\nabla u|^{2}dvol
+B \int_M |u|^{2}dvol.
\ea 
Assume $Ric \ge -a^2 g$ with $a \ge 0$.  Let $1<p<\mu$.  Then there  
holds
\ba \label{SD1-pp}
\|u\|_{\frac{\mu p}{\mu-p}} \le C(\mu, A, B, p) (\|\nabla u\|_{p}+(1+a) \|u\|_p)
\ea
for all $u \in W^{1,p}(M)$, where the constant $C(\mu, A, B, p)$ 
can be bounded from above in terms of upper bounds for $A$, $B$, $\mu$,  
$\frac{1}{\mu-p}$ and $\frac{1}{p-1}$.  } \\

We formulate the corresponding results for the Ricci flow in the situation of [Y1].
The results in the situations of [Y2] and [Y3] can be formulated in a similar way. \\
Consider a smooth solution $g=g(t)$ of the Ricci flow on $M \times [0, T)$, with initial metric 
$g_0$, where 
$M$ is a closed  manifold of dimension $n \ge 3$. Let $\kappa=\kappa(t)$ denote 
$(-\min\{0, \min Ric\})^{1/2}$ at time $t$. \\  

\noindent {\bf Theorem D2} {\it Assume $\lambda_0(g_0)>0$.  Let $1<p<n$.   There is a  positive constant $C$ depending only on 
the dimension $n$, a positive lower bound for $\lambda_0(g_0)$, a positive lower bound for $vol_{g_0}(M)$, an upper bound for $C_S(M,g_0)$, an upper bound for 
$\frac{1}{p-1}$, and an upper bound for $\frac{1}{n-p}$, such that for each $t \in [0, T)$ and 
all $u \in W^{1,p}(M)$ there holds 
\ba 
\|u\|_{\frac{np}{n-p}} \le C (1+R_{max}^+)^{\frac{1}{2}} (\|\nabla u\|_{p}+(1+\kappa)\|u\|_p).
\ea
} \\

\noindent {\bf Theorem D3} {\it Assume $T<\infty$ and $1<p<n$.  There is a positive constant $C$ depending only on the dimension $n$, 
a nonpositive lower bound for $R_{g_0}$, a positive lower bound for $vol_{g_0}(M)$, an upper bound for $C_S(M,g_0)$,
an upper bound for $T$, an upper bound for $\frac{1}{p-1}$, and an upper bound 
for $\frac{1}{n-p}$,  such that for each $t \in [0, T)$ and 
all $u \in W^{1,p}(M)$ there holds 
\ba
\|u\|_{\frac{np}{n-p}} \le C (1+R_{max}^+)^{\frac{1}{2}} (\|\nabla u\|_{p}+(1+\kappa)\|u\|_p).
\ea
} \\

One should compare the above results and the results in Part E with Gallot's estimates of the isoperimetric constant [G1][G2]
which imply estimates for the Sobolev inequalities. In contrast to Gallot's estimates, no upper bound
for the diameter nor positive lower bound for the volume of $g(t)$ is assumed. \\

\noindent {\bf Part E} \\

Based on the $L^p$ estimates of Riesz transforms due to X.~D.~Li [L] we obtain 
a variant of the results in Part F  in the case $1<p<2$. The nonpositive lower bound for the Ricci curvature 
is replaced by the $L^{\frac{n}{2}+\epsilon}$ bound for the (adjusted) negative part of the 
Ricci curvature, where $\epsilon>0$. For a Riemannian manifold $(M, g)$ we set 
$Ric_{min}(x)=\min\{Ric(v,v): v \in T_xM, |v|=1\}$. \\

\noindent {\bf Theorem E1} {\it  Let $(M,g)$ be a complete manifold (without boundary) of dimension 
$n \ge 3$.   Assume  the Sobolev inequality 
\ba \label{SE1}
\left(\int_M |u|^{\frac{2n}{n-2}}dvol\right)^{\frac{n-2}{n}} \le A\int_M |\nabla u|^{2}dvol
+B \int_M |u|^{2}dvol.
\ea 
Let $c\ge 0$ and $\epsilon>0$. 
Assume $(Ric_{min}+c)^- \in L^{\frac{n}{2}+\epsilon}(M)$.    Let $1<p<2$.  Then there  
holds
\ba \label{SE1-pp}
\|u\|_{\frac{n p}{n-p}} \le C (\|\nabla u\|_{p}+(1+\gamma)\|u\|_p)
\ea
for all $u \in W^{1,p}(M)$, where 
\ba
\gamma=\left(\int_M [(Ric_{min}+c)^-]^{\frac{n}{2}+\epsilon}dvol \right)^{\frac{1}{2\epsilon}},
\ea
and the constant $C$ 
can be bounded from above in terms of upper bounds for $n$, $A$, $B$, $c$ and $\frac{1}{\epsilon}$.   
and $\frac{1}{p-1}$.  } \\

\noindent {\bf Theorem E2} {\it Assume $\lambda_0(g_0)>0$.  Let $\epsilon>0$ 
and $1<p<2$. 
There is a  positive constant $C$ depending only on 
the dimension $n$, a positive lower bound for $\lambda_0(g_0)$, a positive lower bound for $vol_{g_0}(M)$, an upper bound for $C_S(M,g_0)$, an upper bound for 
$\frac{1}{p-1}$ and an upper bound for $\frac{1}{\epsilon}$, such that for each $t \in [0, T)$ and 
all $u \in W^{1,p}(M)$ there holds 
\ba
\|u\|_{\frac{np}{n-p}} \le C (1+R_{max}^+)^{\frac{1}{2}} (\|\nabla u\|_{p}+(1+\gamma)\|u\|_p),
\ea
where 
\ba
\gamma=\gamma(t)=\left(\int_M [(Ric_{min}-\frac{1}{n}\min R_{g_0}^-)^-]^{\frac{n}{2}+\epsilon}dvol \right)^{\frac{1}{2\epsilon}}.
\ea
} \\

\noindent {\bf Theorem E3} {\it Assume $T<\infty$. Let  $\epsilon>0$ and $1<p<2$.  There is a positive constant $C$ depending only on the dimension $n$, 
a nonpositive lower bound for $R_{g_0}$, a positive lower bound for $vol_{g_0}(M)$, an upper bound for $C_S(M,g_0)$,
an upper bound for $T$, an upper bound for $\frac{1}{p-1}$, and an upper bound 
for $\frac{1}{\epsilon}$,  such that for each $t \in [0, T)$ and 
all $u \in W^{1,p}(M)$ there holds 
\ba
\|u\|_{\frac{np}{n-p}} \le C (1+R_{max}^+)^{\frac{1}{2}} (\|\nabla u\|_{p}+(1+\gamma)\|u\|_p),
\ea
where $\gamma=\gamma(t)$ is the same quantity as in Theorem E2.
} \\

\sect{From $W^{1,p}$ for lower $p$ to $W^{1,p}$ for higher $p$}

\begin{theo} \label{p-theorem} Consider a Riemannian manifold $(M, g)$ of dimension $n \ge 2$, with or without 
boundary.  Let $1 \le p_0<n$. Assume that the Sobolev inequality 
\ba
\left(\int_M |u|^{\frac{np_0}{n-p_0}}dvol\right)^{\frac{n-p_0}{n}} \le 
A \int_M |\nabla u|^{p_0}dvol 
+\frac{B}{vol_g(M)^{\frac{2}{n}}}\int_M |u|^{p_0}dvol 
\ea
holds true for all $u \in W^{1, p_0}(M)$ with some $A>0$ and $B>0$. Then 
we have 
\ba \label{ppp}
 \left(\int_M |u|^{\frac{pn}{n-p}}dvol \right)^{\frac{n-p}{n}} 
 &\le& 2^{\frac{p-p_0}{p_0}}  A^{\frac{p}{p_0}}(r_p^{p_0}+B)^{\frac{p}{p_0}} \int_M |\nabla u|^p  
 + \frac{2^{\frac{p-p_0}{p_0}}B^{\frac{2p}{p_0}}}{vol_g(M)^{\frac{p}{n}}} \int_M u^p \nonumber \\
 \ea
 for each $p_0<p \le  \frac{n^2p_0}{(n-p_0)^2+np_0}$ and all $u \in W^{1, p}(M)$, where $r_p=\frac{p(n-p_0)}{p_0(n-p)}$
 and the notation of the volume is omitted. 
 \end{theo}
\Pf By scaling invariance we can assume $vol_g(M)=1$.  Consider $u \in C_c^{\infty}(M)$ and set $v=|u|^r$ for $r>1$. Then we have 
\ba
\left(\int_M |u|^{\frac{np_0r}{n-p_0}}\right)^{\frac{n-p_0}{n}} 
 \le A r^{p_0} \int_M |u|^{p_0(r-1)} |\nabla u|^{p_0}
 +B\int_M |u|^{p_0r}. \nonumber \\
\ea  
For a given $p_0<p<n$ we choose $r=r_p$ and hence 
$r_p-1=\frac{n(p-p_0)}{p_0(n-p)}$. Then we infer by H\"{o}lder's inequality
\ba
\left(\int_M |u|^{\frac{np}{n-p}}\right)^{\frac{n-p_0}{n}} 
 &\le&  Ar_p^{p_0} \int_M |u|^{\frac{n(p-p_0)}{n-p}} |\nabla u|^{p_0} + 
 B\int_M |u|^{\frac{p(n-p_0)}{n-p}} \nonumber \\
 &\le& Ar_p^{p_0} \left(\int_M |u|^{\frac{np}{n-p}} \right)^{\frac{p-p_0}{p}}
 \cdot \left(\int_M |\nabla u|^p \right)^{\frac{p_0}{p}}+ 
 B\int_M |u|^{\frac{p(n-p_0)}{n-p}}. \nonumber \\
 \ea
  
 Now we assume $p_0<p\le \frac{n^2p_0}{(n-p_0)^2+np_0}$. Then $\frac{p(n-p_0)}{n-p}\le \frac{np_0}{n-p_0)}$. Since $vol_g(M)=1$ we 
 have by H\"{o}lder's inequality 
 \ba
 \int_M |u|^{\frac{p(n-p_0)}{n-p}}  \le \left(\int_M |u|^{\frac{np_0}{n-p_0}} \right)^{\frac{p(n-p_0)^2}{np_0(n-p)}}
 \ea
 and
 \ba
 \int_M |u|^{\frac{np_0}{n-p_0}}  \le \left(\int_M |u|^{\frac{np}{n-p}}\right)^{\frac{p_0(n-p)}{(n-p_0)p}}. 
 \ea
 We deduce
 \ba
 \left(\int_M |u|^{\frac{pn}{n-p}} \right)^{\frac{(n-p_0)}{n}-\frac{p-p_0}{p}} 
 &\le& Ar_p^{p_0} \left(\int_M |\nabla u|^p \right)^{\frac{p_0}{p}} 
 \nonumber \\ && +B \left(\int_M |u|^{\frac{np_0}{n-p_0}}
 \right)^{\frac{p(n-p_0)^2}{np_0(n-p)}-\frac{(n-p_0)(p-p_0)}{p_0(n-p)}}, 
 \ea
 which leads to
 \ba
 \left(\int_M |u|^{\frac{pn}{n-p}} \right)^{\frac{p_0(n-p)}{np}} 
 &\le& Ar_p^{p_0} \left(\int_M |\nabla u|^p \right)^{\frac{p_0}{p}}+B \left(\int_M |u|^{\frac{np_0}{n-p_0}}
 \right)^{\frac{n-p_0}{n}} \nonumber \\ 
 &\le& Ar_p^{p_0} \left(\int_M |\nabla u|^p \right)^{\frac{p_0}{p}}+AB \int_M |\nabla u|^{p_0} 
 +B^2 \int_M u^{p_0} \nonumber \\
 &\le& A(r_p^{p_0}+B) \left(\int_M |\nabla u|^p \right)^{\frac{p_0}{p}} 
 +B^2 \left(\int_M u^p\right)^{\frac{p_0}{p}}. \nonumber 
 \ea
 It follows that 
 \ba
 \left(\int_M |u|^{\frac{pn}{n-p}} \right)^{\frac{n-p}{n}} 
 &\le& 2^{\frac{p-p_0}{p_0}}  A^{\frac{p}{p_0}}(r_p^{p_0}+B)^{\frac{p}{p_0}} \int_M |\nabla u|^p 
 +2^{\frac{p-p_0}{p_0}}B^{\frac{2p}{p_0}} \int_M u^p. 
 \ea
 By approximation, this holds for all $u \in W^{1,p}(M)$.
 \qed \\
 
 \noindent {\bf Remark} We can also consider the following assumption 
 \ba
\left(\int_M |u|^{\frac{np_0}{n-p_0}}dvol\right)^{\frac{n-p_0}{n}} \le 
A \int_M |\nabla u|^{p_0}dvol 
+\int_M f |u|^{p_0}dvol 
\ea 
for a function $f$. It is easy to adapt the above proof to obtain Sobolev 
inequalitities for higher $p$ in terms of an $L^q$ bound of $f$ for a suitable $q$. 
This can be applied to the Ricci flow to yield Sobolev inequalitites in terms 
of an $L^q$ bound of the scalar curvature.  \\
 
 \begin{lem} \label{p-lemma} Let $1\le p_0<n$. We set $p_{k+1}=\frac{n^2p_k}{(n-p_k)^2+np_k}$ for $k\ge 0$. 
 Then $1\le p_k <n$ for all $k$. Moreover, the sequence $p_k$ is increasing and converges to 
 $n$. 
 \end{lem} 
 \Pf The inequality $p_{k+1} <n$ is equivalent to $(n-p_k)^2>0$, while the 
 inequality $p_k\ge 1$ is equivalent to $(n^2+n)p_k+p_k^2 \ge n^2$. Hence $1\le p_k<n$ follows from 
 the induction. Since $p_k<n$, we have $(n-p_k)^2+np_k<n^2$, and hence $p_{k+1}>p_k$. 
 Let $p_*$ denote the limit of $p_k$. Then $p_*=\frac{n^2 p_*}{(n-p_*)^2+np_*}$. It follows that 
 $p_*=n$. \qed \\

 \noindent {\bf Proof of Theorem A1} Applying Theorem \ref{p-theorem} repeatedly, starting with 
 $p_0=2$.  By induction and 
 Lemma \ref{p-lemma} we then arrive at the desired Sobolev inequalities. 
 \qed \\

 \noindent {\bf Proof of Theorem A2}  We first  observe the following property of the 
 inequality (\ref{pp}): if $A=\alpha A_1$ and $B=\alpha B_1$ for some $\alpha \ge 1$, then 
 we have 
 \ba
 C_1(n, p_0, p, A, B) \le \alpha^{\frac{m(p)p}{p_0}} C_1(n, p_0, p, A_1, B_1)    
 \ea
 and 
 \ba
 C_2(n, p_0, p, A, B) \le \alpha^{\frac{m(p)p}{p_0}} C_1(n, p_0, p, A_1, B_1),
 \ea
 where $m(p)=2^{k+1}$ for $p \in (p_k, p_{k+1}]$ (see Lemma \ref{p-lemma} for $p_k$).
 This follows from the formula (\ref{ppp}). 
 By [Theorem D, Y1], the Sobolev inequality (\ref{sobolev1}) holds true, where $A$ has 
 the dependence as stated in Theorem A2, without reference to $p$.  We then have 
 \ba
 \left(\int_M |u|^{\frac{2n}{n-2}} dvol\right)^{\frac{n-2}{n}} \le A(1+\max R^+ vol(M)^{\frac{2}{n}}) \int_M (|\nabla u|^2 
+\frac{u^2}{vol(M)^{\frac{2}{n}}})dvol.
\ea
Applying Theorem A1 and the above observation we then arrive at the desired Sobolev inequalities. \qed \\

Theorem A3 and Theorem A4 can be proved in the same way.

\sect{Nonlocal Sobolev inequalities in terms of the $(1,p)$-Bessel norm}

First we extend the general results in [Y1] on the heat semigroup and the nonlocal Sobolev inequalitities
to general metrically complete manifolds with or without boundary. 
Consider a Riemannian manifold $(M,g)$ of dimension $n \ge 2$, and 
a function $\Psi \in L^{\infty}(M)$. We set as in [Y1] 
$H=-\Delta+\Psi$ and $Q(u)=\int_M (|\nabla u|^2+\Psi u^2)dvol.$

\begin{theo} \label{D-1} Let $(M, g)$ be metrically complete manifold possibly with boundary. 
Let $0<\sigma^*\le \infty$. Assume that for each $0<\sigma<\sigma^*$ the logarithmic Sobolev 
inequality 
\ba \label{Dlog1}
\int_M u^2 \ln u^2 dvol \le \sigma Q(u)+ \beta(\sigma)
\ea
holds true for all $u \in W^{1,2}(M)$ with $\|u\|_2=1$, where 
$\beta$ is a non-increasing continuous function.  Assume that 
\ba
\tau(t)=\frac{1}{2t}\int^t_0 \beta(\sigma)d\sigma
\ea
is finite for all $0<t < \sigma^*$. Then there holds 
\ba \label{heat1}
\|e^{-tH}u\|_{\infty} \le e^{\tau(t)-\frac{3t}{4}\inf \Psi^-} \|u\|_2
\ea
for each $0<t< \frac{1}{4}\sigma^*$ and all $u \in L^2(M)$. 
There also holds 
\ba \label{heat2}
\|e^{-tH}u\|_{\infty} \le e^{2\tau(\frac{t}{2})-\frac{3t}{4} \inf \Psi^-} \|u\|_1
\ea
for each $0<t< \frac{1}{4}\sigma^*$ and all $u \in L^1(M)$. 
 \end{theo}

\begin{theo}  \label{general} Let $(M,g)$ be a metrically complete manifold, possibly with boundary.
 1) Let $\mu>1$. Assume that $\Psi\ge 0$ and  
for some $c>0$ the inequality 
\ba \label{heatcondition}
\|e^{-tH}u\|_{\infty} \le c t^{-\frac{\mu}{4}} \|u\|_2
\ea 
holds true for each $t>0$ and all $u\in L^2(M)$.  Let $1<p<\mu$.  Then  
there holds 
\ba \label{Hestimate1}
\|H^{-\frac{1}{2}}u\|_{\frac{\mu p}{\mu-p}} \le C(c, \mu, p) \|u\|_p  
\ea
for all $u \in L^p(M)$, where the positive constant 
 $C(\mu, c, p)$ can be bounded from above in terms of upper bounds for $c$, $\mu$,  
$\frac{1}{\mu-p}$ and $\frac{1}{p-1}$.  Consequently, there holds
\ba \label{Hestimate1-2}
\|u\|_{\frac{2\mu}{\mu-2}} \le C(c, 2, p) \left( \int_M (|\nabla u|^2+\Psi u^2)dvol\right)^{\frac{1}{2}} 
\ea
for all $u \in W^{1,2}(M)$. Moreover, there holds for a given $1<p<\infty$
\ba \label{Hestimate1a}
\|u\|_{\frac{\mu p}{\mu-p}} \le C(c, \mu, p) \|H^{\frac{1}{2}}u\|_p
\ea
for all $u \in W^{1,p}(M)$, provided that $(M, g)$ is $p$-Bessel.
\\
2) Let $\mu>1$. Assume that for some $c>0$ the inequality  
\ba 
\|e^{-tH}u\|_{\infty} \le c t^{-\frac{\mu}{4}}\|u\|_2
\ea 
holds true for each $0<t<1$ and all $u\in L^2(M)$. Set $H_0=H-\inf \Psi^-+1$.  Let $1<p<\mu$.  Then there  
holds
\ba \label{Hestimate2}
\|H_0^{-\frac{1}{2}}u\|_{\frac{\mu p}{\mu-p}} \le C(\mu,c, p) \|u\|_p
\ea
for all $u \in L^p(M)$, where the positive constant $C(\mu,c, p)$ has the same property as 
the $C(\mu, c, p)$ above. Consequently, there holds 
\ba \label{Hestimate2-2}
\|u\|_{\frac{2\mu}{\mu -2}} \le C(\mu, 2, p) \left(\int_M (|\nabla u|^2 +(\Psi-\inf \Psi^-+1) u^2)dvol\right)^{\frac{1}{2}}
\ea
for all $u \in W^{1,2}(M)$. Moreover, there holds for a given $1<p<\infty$
\ba \label{Hestimate2a}
\|u\|_{\frac{\mu p}{\mu -p}} \le C(\mu, c, p) \|H_0^{\frac{1}{2}}u\|_p
\ea
for all $u \in W^{1,p}(M)$, provided that $(M, g)$ is $p$-Bessel.
\end{theo}

To establish these two results, we need the following two ingredients: the construction of the 
heat semigroup $e^{-tH}$ and the $L^p$ contraction properties of $e^{-tH}$ for all $1<p<\infty$. 
In [Y1], since the manifold is assumed to be closed, the heat semigroup is constructed by using the spectral representation in terms of the eigenfunctions. This works equally well on a compact manifold with 
boundary, where the eigenfunctions satisfy the Neumann boundary condition.  For a general 
metrically complete manifold, we follow the construction in [St] based on the general theory of 
spectral representation of self-adjoint operators.  The case of $H=-\Delta$ on a complete 
manifold without boundary is treated in [St], but the arguments extend quite easily to general 
metrically complete manifolds and $H=-\Delta+\Psi$ with $\Psi\in L^{\infty}(M)$ and $\Psi \ge 0$.

Consider $\Psi \in L^{\infty}(M)$ with $\Psi \ge 0$.  The initial domain for $H=-\Delta+\Psi$ is the space $\Omega_H=C^{\infty}_{c,N}(M)=\{
u \in C^{\infty}_c(M): \frac{\partial u}{\partial \nu}=0\}$, where $\nu$ denotes the 
inward unit normal of $\partial M$, which is dense in $L^2{M}$.  Let $H_{min}$ 
denote the $L^2$ closure of $H$, whose domain $D(H_{min})$ consists of all $u \in L^2(M)$ such that 
there is a sequence $u_i \in \Omega_H$ such that $u_i \rightarrow u$ in $L^2(M)$ and 
$Hu_i$ converges in $L^2(M)$ to some function, which we can write $Hu$.    
Let $H_{max}$ be the adjoint of $H_{min}$ in $L^2(M)$, and $D(H_{max}) \subset L^2(M)$ its domain. 
 We have the following extension 
of [St, Lemma 2.3]. 

\begin{lem} \label{0lemma} Let $(M, g)$ be metrically complete. Assume that $u \in D(H_{max})$ 
satisfies $Hu=\lambda u$ for some $\lambda<0$.   Then $u \equiv 0$.
\end{lem}      
\Pf   By basic elliptic regularity we have $u\in  W^{2,p}_{loc}(M)$ for 
all $p>0$ and $\frac{\partial u}{\partial \nu}=0$.  ( By the Sobolev embedding we have 
$u \in C^1(M)$. ) Fix $x_0 \in M$ and let $\varphi(x)=\varphi_{r_1, r_2}(x)
=\psi((r_2-r_1)^{-1}(d(x_0,x)+r_2-2r_1))$ for a smooth function 
$\psi(t)$  which is $1$ for $t\le 1$ and $0$ for $t\ge 2$.  Then 
$|\nabla \varphi| \le c(r_2-r_1)^{-1}$ for a constant $c$.  

Now we have 
\ba
\lambda <\varphi^2u, u>_2&=& <\varphi^2 u, Hu>_2 = -<\varphi^2 u, \Delta u>_2
+<\varphi^2 u, \Psi u> \nonumber \\
&\ge& -<\varphi^2 u, \Delta u>_2 =\|\varphi^2\nabla u\|_2^2+2<u \nabla \phi, \phi \nabla u>_2.  
\ea 
It follows that  
\ba 
\|\varphi^2\nabla u\|_2^2 \le \lambda <\varphi^2 u, u>_2 
+2<u \nabla \phi, \phi \nabla u>_2
\ea
By Schwarz inequality we then deduce
\ba
\|\varphi^2\nabla u\|_2^2 \le 2\lambda <\varphi^2 u, u>_2 
+\frac{4c^2}{(r_2-r_1)^2}\|u\|^2_2.
\ea
Letting first $r_2 \rightarrow \infty$ and then $r_1 \rightarrow \infty$ we arrive at 
\ba
\|\nabla u\|^2_2 \le 2\lambda \|u\|_2^2.
\ea 
Since $\lambda<0$, we conclude $u \equiv 0$. \qed \\

By this lemma and [St, Lemma 2.1] we infer that $H_{max}=H_{min}$, which is 
the self-adjoint extension of $H$. Now we can apply the spectral theorem for 
self-adjoint operators to obtain the heat semigroup $e^{-tH}$ and other 
potentials of $H$ such as $H^{-\frac{1}{2}}$ and $H^{\frac{1}{2}}$.    

In [Y1], the $L^p$ contraction property of $e^{-tH}$ is derived in terms of the 
$L^2$ contraction property and the $L^{\infty}$ contraction property, with the latter 
implied by the maximum principle. This argument can be extended to compact manifolds 
with boundary. But the maximum principle may not hold on noncompact manifolds. 
Instead, we follow the arguments in [St] for obtaining the $L^p$ contraction property.  By 
the arguments in Section 3 of [St], 
in order to show that $e^{-tH}$ is a contraction on $L^p(M) \cap L^2(M)$ for each $1 \le p \le \infty$, it suffices to 
establish the following two lemmas. 

\begin{lem}   For each $1<p<\infty$ the operator $H$ with domain $\Omega_H$ is dissipative, i.e. 
for each nonzero $u \in \Omega_H$, there is a function $v \in L^q$ with $q=\frac{p}{p-1}$ such that 
$\|v\|_q=\|u\|_p, <u,v>_2=\|u\|_p$ and $<Hu, v>_2 \le 0$. 
\end{lem} 

\begin{lem}
Let $1<p\le q<\infty$. Assume that $u \in L^p(M) \cap L^q(M)$ satisfies
$Hu=\lambda u$ for some $\lambda<0$ (this contains the assumption that 
$u$ lies in the domain of the closure of $H$ in $L^p(M)$ and that in $L^q(M)$.)  Then $u\equiv 0$.
\end{lem}

Since $\Psi\ge 0$, the proof of Lemma 3.1 and the proof of Lemma 3.4 carry over.  The treatment of 
$\Psi$ here is similar to  
to the above proof of Lemma \ref{0lemma}.

Having established the desired construction of $e^{-tH}$ and the $L^p$ contraction properties we make 
two more remarks. First, the construction of the heat semigroup $e^{-tH}$ for a general 
$\Psi \in L^{\infty}(M)$ follows via the formula $e^{-tH}=e^{-t \inf \Psi^-} e^{-tH_1}$, 
where $H_1=-\Delta+\Psi-\inf \Psi^-$. Second, in [Y1] the space $L^{\infty}(M)$ is used 
in the formulations of the Marcinkiewicz 
interpolation theorem and the Riesz-Thorin interpolation theorem.  In the case of 
a general Riemannian manifold we replace $L^{\infty}(M)$ by 
$L^{\infty}(M) \cap L^1(M)$.  \\

\noindent {\bf Proof of Theorem \ref{D-1}} Consider $u_0 \in L^2(M)$. We claim that 
$e^{-tH}u_0 \in W^{1,2}(M)$ for $t>0$. Indeed we have for $u=e^{-tH}u_0$ 
\ba
\frac{\partial u}{\partial t}=Hu, \, \, \, \frac{\partial u}{\partial \nu}=0.
\ea
Then we have 
\ba
\frac{d}{dt} \int_M \varphi^2 u^2 =2\int_M \varphi^2 u Hu =
-2\int_M \varphi^2 |\nabla u|^2 -2\int_M \varphi u \nabla u \cdot \nabla \varphi -2\int_M \Psi \varphi^2 
u^2, \nonumber \\
\ea
where $\varphi=\varphi_{r_2,r_1}$ is the function in the proof of Lemma \ref{0lemma}. It follows that
\ba
\frac{d}{dt} \int_M \varphi^2 u^2 \le 
-\int_M \varphi^2 |\nabla u|^2 +\int_M u^2 |\nabla \varphi|^2 -\int_M \Psi \varphi^2 
u^2, 
\ea
and then
\ba
\int_M \varphi^2 u^2 +\int_0^t \int_M \varphi^2 |\nabla u|^2 \le 
\int_M \varphi^2 u^2 |_{t=0} + \int_0^t \int_M (|\nabla \varphi|^2-\Psi)u^2. 
\ea  
Letting $r_2 \rightarrow \infty$ and then $r_1 \rightarrow \infty$ we arrive at 
\ba
\int_M  u^2 +\int_0^t \int_M |\nabla u|^2 \le 
\int_M  u^2 |_{t=0} -\int_0^t \int_M \Psi u^2. 
\ea
It follows that $\int_M |\nabla u(\cdot, t)|^2 <\infty$ for a.e. $t>0$. By continuity, 
$\int_M |\nabla u(\cdot, t)|^2<\infty$ for all $t>0$. Hence $e^{-tH}u_0 \in W^{1,2}(M)$ 
for all $t>0$. 

Now we can carry over the proof of Theorem 5.3 in [Y1]. Some modification is necessary 
because $M$ is possibly noncompact. Let $\varphi=\varphi_{r_1, r_2}$ be the function in the 
proof of Lemma \ref{0lemma}. In place of [Y1, (B.13)] we have now for $u_s=e^{-sH}u_0$ for 
a given $u_0 \in W^{1,2}(M) \cap L^{\infty}(M)$ 
\be
\frac{d}{ds}\ln (e^{-N(s)}\|\varphi u_s \|_{p(s)})=\frac{d}{ds}\left(-N(s)+\frac{1}{p(s)} 
\ln \|\varphi u_s\|_{p(s)}^{p(s)}\right) \nonumber
\ee
\be
=\frac{\Gamma}{\sigma}-\frac{1}{p^2}\frac{p}{\sigma} \ln \|\varphi u_s\|_p^p +\frac{1}{p} \|\varphi u_s\|_p^{-p}
\left(-pQ(\varphi u_s, u_s^{p-1})+\frac{p}{\sigma} \int_M \varphi u_s^p \ln u_s\right) \nonumber 
\ee
\be
=\frac{1}{\sigma}\|\varphi u_s\|^{-p}_p \left( \int_M \varphi u_s^p \ln u_s -
\sigma Q(\varphi u_s, u_s^{p-1})-\Gamma \|\varphi u_s\|_p^p-\|\varphi u_s\|_p^p \ln \|\varphi u_s\|_p \right).
\ee
It follows that 
\ba \label{monotone}
&& \ln (e^{-N(t_2)}\|\varphi u_{t_2} \|_{p(t_2)}) \le \ln (e^{-N(t_1)}\|\varphi u_{t_1} \|_{p(t_1)}) +
\nonumber \\ 
&& \int_{t_1}^{t_2}  \frac{1}{\sigma}\|\varphi u_s\|^{-p}_p \left( \int_M \varphi u_s^p \ln u_s -
\sigma Q(\varphi u_s, u_s^{p-1})-\Gamma \|\varphi u_s\|_p^p-\|\varphi u_s\|_p^p \ln \|\varphi u_s\|_p \right)ds
\nonumber \\
\ea
for $t_2>t_1>0$. Letting first $r_2 \rightarrow \infty$ and then $r_1 \rightarrow \infty$ we then 
arrive at (\ref{monotone}) without the presence of $\phi$.  
\qed \\

\noindent {\bf Proof of Theorem \ref{general}} Given the $L^p$ contraction property established 
above, the proof of Theorem C.5 in [Y1] carries over 
straightforwardly. The Sobolev inequality (\ref{Hestimate1-2}) follows because of the 
identity $<H^{\frac{1}{2}}u, H^{\frac{1}{2}}u>_2 =\int_M (|\nabla u|^2+\Psi u^2)dvol$ for all 
$u \in C^{\infty}_c(M)$ (then also for $u \in W^{1,2}(M)$).  The Sobolev inequality 
(\ref{Hestimate2-2}) follows in the same fashsion.  \qed \\

\noindent {\bf Proof of Theorem B4} Assume $\int_M u^2 =1$. As in the proof of Theorem 3.1 in [Y1] we have 
\ba
\ln \int_M |u|^{\frac{2\mu}{2-\mu}} =\ln \int_M  u^2 |u|^{\frac{4}{\mu-2}} \ge \int_M u^2 \ln |u|^{\frac{4}{\mu-2}}.
\ea
It follows that
\ba
\int_M u^2 \ln u^2  &\le& \frac{\mu}{2} \ln \left( \int_M |u|^{\frac{2\mu}{\mu-2}}\right)^{\frac{\mu-2}{\mu}} 
\le \frac{\mu}{2} \ln \left(A\int_M (|\nabla u|^2 +\Psi \int_M u^2 \right) 
\nonumber \\
&\le& \frac{\mu}{2} \ln A +\frac{\mu}{2} \ln \int_M (|\nabla u|^2+\Psi u^2).
 \ea 
By [Y1, Lemma 3.2] we then deduce each $\sigma>0$
\ba
\int_M u^2 \ln u^2 \le  \frac{\mu}{2} \sigma \int_M (|\nabla u|^2+\Psi u^2)-\frac{\mu}{2} \ln \sigma 
+\frac{\mu}{2} \ln A-1, 
\ea
which leads to 
\ba
\int_M u^2 \ln u^2 \le \sigma \int_M (|\nabla u|^2+\Psi u^2)-\frac{\mu}{2} \ln \sigma 
+\frac{\mu}{2} \ln \frac{\mu}{2}+\frac{\mu}{2} \ln A-1.
\ea
By Theorem \ref{D-1} we deduce for $H=-\Delta+1$
\ba
\|e^{-tH}u\|_{\infty} \le  t^{-\frac{\mu}{4}} e^{\frac{\mu}{4}+\frac{A_0}{2}} \|u\|_2
\ea
for all $t>0$, where $A_0=\frac{\mu}{2} \ln \frac{\mu}{2}+\frac{\mu}{2} \ln A-1.$
Applying Theorem \ref{general} we then arrive at the desired inequality (\ref{SB4-p}). 
\qed \\

\noindent {\bf Proof of Theorem B1} The Sobolev inequality (\ref{SB1}) leads to 
\ba
\left(\int_M |u|^{\frac{2\mu}{\mu-2}}dvol \right)^{\frac{\mu-2}{\mu}} \le 
\max\{A, B\} \int_M (|\nabla u\|^2+u^2)dvol.
\ea
Hence we can apply Theorem B4. \qed \\

\begin{lem} \label{changelemma} Let $\bar g=\lambda^2 g$ for some $\lambda \ge 1$. 
Let $\mu>1$ and $1 \le p<\mu$.  Assume the inequality 
\ba
\|u\|_{\frac{\mu p}{\mu-p}} \le C \|u\|_{B, 1, p}
\ea
for all $u \in W^{1, p}(M)$ with respect to $\bar g$. Then there holds 
\ba
\|u\|_{\frac{\mu p}{\mu-p}} \le \lambda C \|u\|_{B, 1, p}
\ea 
for all $u \in W^{1,p}(M)$ with respect to $g$.
\end{lem} 
\Pf We compute the scaling change of $(-\Delta+1)^{\frac{1}{2}}$. 
 We have 
$\Delta_{\bar g}=\lambda^{-2} \Delta$. Hence 
\ba
-\Delta_{\bar g}+1=-\lambda^{-2} \Delta+1=\lambda^{-2}(-\Delta+\lambda^2).
\ea
By [B, Lemma 4.2],  we have for $u \in L^p(M)$ and $a\ge 0$
\ba \label{equivalence}
c_1(a \|u\|_p+\|(-\Delta)^{\frac{1}{2}}u\|_p) \le \|(-\Delta+a^2)^{\frac{1}{2}}u\|_p 
\le c_2(a \|u\|_p + \|(-\Delta)^{\frac{1}{2}}u\|_p),
\ea
where $c_1$ and $c_2$ are universal constants.
It follows that 
\ba
\|(-\Delta+\lambda^2)^{\frac{1}{2}} u\|_p 
&\le& c_2(\lambda\|u\|_p+\|(-\Delta)^{\frac{1}{2}}u\|_p) 
\nonumber \\
&\le& c_2\lambda(\|u\|_p+\|(-\Delta)^{\frac{1}{2}}u\|_p) 
\nonumber \\
&\le& c_1c_2 \lambda \|(-\Delta+1)^{\frac{1}{2}}u\|_p.
\ea
Hence 
\ba
\|(-\Delta_{\bar g}+1)^{\frac{1}{2}}u\|_p \le c_1c_2 \|(-\Delta+1)^{\frac{1}{2}}u\|_p.
\ea
Now we have 
\ba
\|u\|_{\frac{pn}{n-p}, \bar g}=\|u\|_{\frac{np}{n-p}} \lambda^{\frac{n-p}{p}}
\ea
and 
\ba
\|(-\Delta_{\bar g}+1)^{\frac{1}{2}}u\|_{p, \bar g} = \|(-\Delta_{\bar g}+1)^{\frac{1}{2}}u\|_p \lambda^{\frac{n}{p}}.
\ea 
We arrive at 
\ba
\|u\|_{\frac{np}{n-p}} \le \lambda C\|(-\Delta+1)^{\frac{1}{2}}u\|_p.
\ea
\qed \\

\noindent {\bf Proof of Theorem B2} By [Y1, Theorem D], the Sobolev inequality 
(\ref{sobolev1}) holds true, where $A$ has the same property as the $C$ in the theorem without 
the reference to $p$. Set $\lambda=\lambda(t)=(1+R_{max}^+)^{1/2}$ at time $t$.  Then the Sobolev 
inequality (\ref{sobolev1}) still holds true for $\bar g=\lambda^2 g$. Since $R^+_{max} 
\le 1$ for $\bar g$, we deduce  
\ba
\left(\int_M |u|^{\frac{2n}{n-2}}dvol\right)^{\frac{n-2}{n}} \le
A\int_M (|\nabla u|^2+ u^2)dvol.
\ea
Applying Theorem B1 and Lemma  \ref{changelemma} we then arrive at the desired inequality 
(\ref{SB2-p}). \qed \\

\noindent {\bf Proof of Theorem B3} By [Y1, Theorem $\mbox{D}^*$], the Sobolev 
inequality (\ref{sobolev2}) holds true, where $A$ and $B$ have the same property as 
the $C$ in the theorem, without the reference to $p$. Let $\lambda$ and $\bar g$ be the same 
as above. Then we have for $\bar g$
\ba
\left(\int_M |u|^{\frac{2n}{n-2}}dvol\right)^{\frac{n-2}{n}} &\le&
A\int_M (|\nabla u|^2+\frac{R}{4}u^2) dvol+\frac{B}{\lambda^2} \int_M u^2dvol 
\nonumber \\ 
&\le& A \int_M |\nabla u|^2 dvol+ (\frac{A}{4}+B)\int_M u^2 dvol \nonumber \\
&\le& (A+B) \int_M (|\nabla u|^2+u^2)dvol.
\ea
Applying Theorem B1 and Lemma \ref{changelemma} we arrive at the desired 
inequality (\ref{SB3-p}). \qed \\

\sect{$W^{2,p}$ Sobolev inequalities}

\noindent {\bf Proof of Theorem C1}  Let $u \in W^{2,p}(M)$ for $1<p<\frac{\mu}{2}$. 
Since $(-\Delta+\Psi)^{\frac{1}{2}}$ is a pseudo-differential operator of order $1$ 
[Se] on a compact manifold, it is a bounded map from $W^{2,p}(M)$ into $W^{1,p}(M)$. 
Hence $v=(-\Delta+\Psi)^{\frac{1}{2}}u \in W^{1,p}(M)$. Applying Theorem B4 to
$v$ we infer 
\ba
\|v\|_{\frac{\mu p}{\mu-p}} \le C(\mu, A, p) \|(-\Delta+\Psi)^{\frac{1}{2}}v\|_p =C(\mu, A, p) \|(-\Delta+\Psi)u\|_p,
\ea
i.e.
\ba
\|(-\Delta+\Psi)^{\frac{1}{2}}u\|_{\frac{\mu p}{\mu-p}} \le C(\mu, A, p) \|(-\Delta+\Psi)u\|_p.
\ea
For each $1<q<\infty$, $(-\Delta+\Psi)^{\frac{1}{2}}$ is a bounded operator from 
$W^{1, q}(M)$ into $L^q(M)$ with the bounded inverse $(-\Delta+\Psi)^{-\frac{1}{2}}$.  
Hence we deduce $u \in W^{1, \frac{\mu p}{\mu-p}}(M)$. Applying Theoem E4 to $u$ with the exponent 
$\frac{\mu p}{\mu -p}$ instead of $p$ we then infer 
\ba
\|u\|_{\frac{\mu p}{\mu -2p}} &\le& C(\mu, A, \frac{\mu p}{\mu -p}) 
\|(-\Delta+\Psi)^{\frac{1}{2}}u\|_{\frac{\mu p}{\mu-p}} 
\nonumber \\ &\le& C(\mu, A, \frac{\mu p}{\mu -p}) 
C(\mu, A, p) \|(-\Delta+\Psi)u\|_p.
\ea
(Note that $1<\frac{\mu p}{\mu-p}<\mu$ because $1<p<\frac{\mu}{2}$.)
\qed \\

Theorem C2 and Theorem C3 follow from Theorem C1, and [Y1, Theorem D] and 
[Y1, Theorem $\mbox{D}^*$] respectively.

\sect{Estimates of the Riesz transform and $W^{1,p}$ Sobolev inequalities} 

The following theorem is a consequence of D.~Bakry's result on 
$L^p$ estimates for the Riesz transform [B]. 

\begin{theo} \label{bk-estimate}Let $(M, g)$ be a complete Riemannian manifold (without boundary) of 
dimension $n \ge 2$ such that the Ricci curvature is bounded from below by 
$-a^2$ for some $0\le a <\infty$.  Then there holds for each $1<p<\infty$
\ba \label{bk-p}
\|(-\Delta+1)^{\frac{1}{2}}u\|_p \le C(p)(\|\nabla u\|_p+(1+a)\|u\|_p)
\ea
for all $u \in W^{1,p}(M)$, where the constant $C(p)$ depends only on $p$.  
\end{theo} 
\Pf In [B] the operator $-\Delta+\nabla \phi \cdot \nabla$ for a given function $\phi$ is handled. 
It is easy to see that all the arguments in [B] go through for the operator 
$-\Delta+1$.   Hence [B, Theorem 4.1] extends to yield for $1<q<\infty$
\ba \label{bk}
\|\nabla v\|_q \le C_q (\|(-\Delta+1)^{\frac{1}{2}}v\|_q+a\|v\|_q)
\ea
for all $u \in C^{\infty}_c(M)$, where $C_q$ depends only on $q$. 
On the other hand, we have  $\|e^{t(\Delta-1)}v\|_q \le e^{-t}\|v\|_q$ for $1<q<\infty$ and all $v\in L^q(M) \cap L^2(M)$.  Applying  
the formula $(-\Delta+1)^{-\frac{1}{2}}=\Gamma(\frac{1}{2})^{-1}\int_0^{\infty} t^{-\frac{1}{2}}e^{t(\Delta-1)}dt$
we infer $\|(-\Delta+1)^{-\frac{1}{2}}v\|_q \le \|v\|_q$ and hence 
\ba \label{q-con}
\|v\|_q \le \|(-\Delta+1)^{\frac{1}{2}}v\|_q
\ea
for all $v \in L^q(M)\cap L^2(M)$.  

 Since $(-\Delta+1)^{\frac{1}{2}}(C^{\infty}_c(M))$ is dense 
in $L^q(M)$ (see [CD][R]), we have for $u \in C^{\infty}_c(M)$, $1<p<\infty$ and $q=\frac{p}{p-1}$
\ba \label{dual-1}
\|(-\Delta+1)^{\frac{1}{2}}u\|_p &=& \sup \{<(-\Delta +1)^{\frac{1}{2}}u, 
(-\Delta+1)^{\frac{1}{2}}v>_2: \nonumber \\
&& v \in C^{\infty}_c(M), \|(-\Delta+1)^{\frac{1}{2}}v\|_q \le 1\}.
\ea 
But 
\ba \label{dual-2}
<(-\Delta +1)^{\frac{1}{2}}u, (-\Delta+1)^{\frac{1}{2}}v>_2&=&<(-\Delta+1)u, v>_2=
\int_M \nabla u \cdot \nabla v + \int_M uv \nonumber \\
&\le& \|\nabla u\|_p \|\nabla v\|_q+ \|u\|_p \|v\|_q  \nonumber \\
&\le& (\|\nabla u\|_p+\|u\|_p) (\|\nabla v\|_q+\|v\|_q).
\ea
By (\ref{q-con}) and (\ref{bk}) we then deduce $\|\nabla v\|_q+\|v\|_q 
\le (1+(a+1)C(q))\|(-\Delta+1)^{\frac{1}{2}}v\|_q$.  
By (\ref{dual-1}) and (\ref{dual-2}) we then arrive at 
\ba \label{p-estimate}
\|(-\Delta+1)^{\frac{1}{2}}u\|_p \le  C(p,a)(\|\nabla u\|_p+\|u\|_p),
\ea
where $C(p,a)=1+(a+1)C_{\frac{p}{p-1}}$.
By (\ref{bk}, (\ref{q-con}) (applied to $p$) and (\ref{p-estimate}) we conclude that 
$(M, g)$ is $(1,p)$-Bessel and that (\ref{p-estimate}) holds true for all $u \in W^{1,p}(M)$. 

To derive the inequality (\ref{bk-p}), we consider the metric $\bar g=\lambda^2 g$, where 
$\lambda=1+a$. Since the Ricci curvature of $\bar g$ is bounded from below by $-\frac{a^2}{(1+a)^2}
\ge -1$, we have by (\ref{p-estimate})
\ba
\|(-\Delta_{\bar g}+1)^{\frac{1}{2}}u\|_{p, \bar g} \le C(p, 1)(\|\nabla_{\bar g} u\|_{p, \bar g}+\|u\|_{p, \bar g})
\ea
for $1<p<\infty$ and all $u \in L^p(M)$.  But $\Delta_{\bar g}=
\lambda^{-2} \Delta$. Hence we obtain 
\ba
\|(-\Delta+\lambda^2)^{\frac{1}{2}}u\|_{p, \bar g} \le \lambda C(p, 1) (\|\nabla_{\bar g} u\|_{p, \bar g}+\|u\|_{p, \bar g}).
\ea
Transforming to $g$ we obtain 
\ba
\|(-\Delta+\lambda^2)^{\frac{1}{2}}u\|_{p} \le C(p, 1) (\|\nabla u\|_{p}+\lambda \|u\|_{p}).
\ea
By (\ref{equivalence}) there holds
\ba
\|(-\Delta+\lambda^2)^{\frac{1}{2}}u\|_{p} &\ge& c_1 (\lambda\|u\|_p
+\|(-\Delta)^{\frac{1}{2}}u\|_p) \ge c_1(\|u\|_p+\|(-\Delta)^{\frac{1}{2}}u\|_p)
\nonumber \\ &\ge& c_1c_2^{-1} \|(-\Delta+1)^{\frac{1}{2}}u\|_p.
\ea
It follows that 
\ba
\|(-\Delta+1)^{\frac{1}{2}}u\|_p\le c_1^{-1}c_2 C(p, 1) (\|\nabla u\|_{p}+\lambda \|u\|_{p})
\ea
which gives rise to (\ref{bk-p}). 
\qed \\

The next result is a consequence of the $L^p$ estimates for the Riesz transform due to 
X.~D.~Li [L]. \\

\begin{theo} \label{L-result} Let $(M,g)$ be a complete Riemannian manifold (without boundary) 
of dimension $n\ge 2$. Assume that there is a constant $c_1$ such that
\ba
\|e^{t(\Delta-1)}u\|_{\infty} \le c_1t^{-\frac{n}{2}}\|u\|_1
\ea
for all $u \in L^1(M)$ and $0<t\le 1$. Assume that $(Ric_{min}+c_2)^- \in L^{\frac{n}{2}+\epsilon}(M)$ for some $c_2\ge 0$ and $\epsilon>0$.  Then there holds for each $1<p<2$
\ba \label{bk-p}
\|(-\Delta+1)^{\frac{1}{2}}u\|_p \le C(\|\nabla u\|_p+(1+\gamma) \|u\|_p)
\ea
for all $u \in W^{1,p}(M)$, where
\ba
\gamma=\left( \int_M [(Ric_{min}+c_2)^-]^{\frac{n}{2}+\epsilon}dvol\right)^{\frac{1}{2\epsilon}},
\ea
 and the constant $C$ can be bounded above in terms of 
upper bounds for $n, c_1, c_2, \frac{1}{\epsilon}$ and $ \frac{1}{p-1}$. 
\end{theo} 
\Pf This follows from the proof of Theorem 2.2 in [L] and the arguments in the above proof of 
Theorem \ref{bk-estimate}. \qed \\

\noindent {\bf Proof of Theorems D1, D2 and D3}  Combine Theorem \ref{bk-estimate} with Theorems B1, B2 and 
B3.  \qed \\

\noindent {\bf Proof of Theorems E1, E2 and E3} Combine Theorem \ref{L-result} with 
Theorems B1, B2 and B3, and also apply Theorem 3.1 and the arguments in the proof of Theorem B4. \qed

\end{document}